\documentclass[a4paper,12pt]{article}

 \usepackage{latexsym,amssymb,amsmath}

 \newtheorem{theorem}{Theorem}

\begin{document}

\title{A Cauchy-Kowalevski theorem for inframonogenic functions\thanks{accepted for publication in Mathematical Journal of Okayama University}}

\author{Helmuth R. Malonek$^{\star,1}$, Dixan Pe\~na Pe\~na$^{\star,2}$\\and Frank Sommen$^{\dagger,3}$}

\date{\normalsize{$^\star$Department of Mathematics, Aveiro University,\\3810-193 Aveiro, Portugal\\
$^\dagger$Department of Mathematical Analysis, Ghent University,\\9000 Gent, Belgium}\\\vspace{0.4cm}
\small{$^1$e-mail: hrmalon@ua.pt\\
$^2$e-mail: dixanpena@ua.pt; dixanpena@gmail.com\\
$^3$e-mail: fs@cage.ugent.be}}

\maketitle

\begin{abstract}
\noindent
In this paper we prove a Cauchy-Kowalevski theorem for the functions satisfying the system $\partial_xf\partial_x=0$ (called inframonogenic functions).\vspace{0.2cm}\vspace{0.1cm}
\end{abstract}

\section{Introduction}

Let $\mathbb{R}_{0,m}$ be the $2^m$-dimensional real Clifford algebra constructed over the orthonormal basis $(e_1,\ldots,e_m)$ of the Euclidean space $\mathbb R^m$ (see \cite{Cl}). The multiplication in $\mathbb{R}_{0,m}$ is determined by the relations $e_je_k+e_ke_j=-2\delta_{jk}$, $j,k=1,\dots,m$, where $\delta_{jk}$ is the Kronecker delta. A general element of $\mathbb{R}_{0,m}$ is of the form
\[a=\sum_Aa_Ae_A,\quad a_A\in\mathbb R,\]
where for $A=\{j_1,\dots,j_k\}\subset\{1,\dots,m\}$, $j_1<\cdots<j_k$, $e_A=e_{j_1}\dots e_{j_k}$. For the empty set $\emptyset$, we put $e_{\emptyset}=1$, the latter being the identity element.

Notice that any $a\in\mathbb{R}_{0,m}$ may also be written as  $a=\sum_{k=0}^m[a]_k$ where $[a]_k$ is the projection of $a$ on $\mathbb R_{0,m}^{(k)}$. Here $\mathbb R_{0,m}^{(k)}$ denotes the subspace of $k$-vectors defined by
\[\mathbb R_{0,m}^{(k)}=\biggl\lbrace a\in\mathbb R_{0,m}:\;a=\sum_{\vert A\vert=k}a_Ae_A,\quad a_A\in\mathbb R\biggr\rbrace.\]
Observe that $\mathbb R^{m+1}$ may be naturally identified with $\mathbb R_{0,m}^{(0)}\oplus\mathbb R_{0,m}^{(1)}$ by associating to any element $(x_0,x_1,\ldots,x_m)\in\mathbb R^{m+1}$ the \lq\lq paravector" $x=x_0+\underline x=x_0+\sum_{j=1}^mx_je_j$.

Conjugation in $\mathbb R_{0,m}$ is given by
\[\overline a=\sum_Aa_A\overline e_A,\]
where $\overline e_A=\overline e_{j_k}\dots\overline e_{j_1}$, $\overline e_j=-e_j$, $j=1,\dots,m$. One easily checks that $\overline{ab}=\overline b\overline a$ for any $a,b\in\mathbb R_{0,m}$. Moreover, by means of the conjugation a norm $\vert a\vert$ may be defined for each $a\in\mathbb R_{0,m}$ by putting
\[\vert a\vert^2=[a\overline a]_0=\sum_Aa_A^2.\]
Let us denote by  $\partial_x=\partial_{x_0}+\partial_{\underline x}=\partial_{x_0}+\sum_{j=1}^me_j\partial_{x_j}$ the generalized Cauchy-Riemann operator and let $\Omega$ be an open set of $\mathbb R^{m+1}$. According to \cite{MDS}, an $\mathbb R_{0,m}$-valued function $f\in C^2(\Omega)$ is called an inframonogenic function in $\Omega$ if and only if it fulfills in $\Omega$ the \lq\lq sandwich" equation $\partial_xf\partial_x=0$.

It is obvious that monogenic functions (i.e. null-solutions of $\partial_x$) are inframonogenic. At this point it is worth remarking that the monogenic functions are the central object of study in Clifford analysis (see \cite{BDS,CnMa,DSS,GuSp,K,KS,Ma,R}). Furthermore, the concept of monogenicity of a function may be seen as the higher dimensional counterpart of holomorphy in the complex plane.

Moreover, as
\[\Delta_x=\sum_{j=0}^m\partial_{x_j}^2=\partial_x\overline\partial_x=\overline\partial_x\partial_x,\]
every inframonogenic function $f\in C^4(\Omega)$ satisfies in $\Omega$ the biharmonic equation $\Delta_x^2f=0$ (see e.g. \cite{BG,GK,M,So}).

This paper is intended to study the following Cauchy-type problem for the inframonogenic functions. Given the functions $A_0(\underline x)$ and $A_1(\underline x)$ analytic in an open and connected set $\underline\Omega\subset\mathbb R^m$, find a function $F(x)$ inframonogenic in some open neighbourhood $\widetilde\Omega$ of $\underline\Omega$ in $\mathbb R^{m+1}$ which satisfies
\begin{align}
F(x)\vert_{x_0=0}&=A_0(\underline x),\label{cond1}\\
\partial_{x_0}F(x)\vert_{x_0=0}&=A_1(\underline x).\label{cond2}
\end{align}

\section{Cauchy-type problem for inframonogenic functions}

Consider the formal series
\begin{equation}\label{series}
F(x)=\sum_{n=0}^\infty x_0^nA_n(\underline x).
\end{equation}
It is clear that $F$ satisfies conditions (\ref{cond1}) and (\ref{cond2}). We also see at once that
\[\partial_x\left(x_0^nA_n\right)\partial_x=n(n-1)x_0^{n-2}A_n+nx_0^{n-1}\big(\partial_{\underline x}A_n+A_n\partial_{\underline x}\big)+x_0^n\partial_{\underline x}A_n\partial_{\underline x}.\]
We thus get
\[\partial_xF\partial_x=\sum_{n=0}^\infty x_0^n\Big((n+2)(n+1)A_{n+2}+(n+1)\big(\partial_{\underline x}A_{n+1}+A_{n+1}\partial_{\underline x}\big)+\partial_{\underline x}A_n\partial_{\underline x}\Big).\]
From the above it follows that $F$ is inframonogenic if and only if the functions $A_n$ satisfy the recurrence relation
\[A_{n+2}=-\frac{1}{(n+2)(n+1)}\Big((n+1)\big(\partial_{\underline x}A_{n+1}+A_{n+1}\partial_{\underline x}\big)+\partial_{\underline x}A_n\partial_{\underline x}\Big),\quad n\ge0.\]
It may be easily proved by induction that
\begin{equation}\label{sol}
A_n=\frac{(-1)^{n+1}}{n!}\left(\sum_{j=0}^{n-2}\partial_{\underline x}^{n-j-1}A_0\partial_{\underline x}^{j+1}+\sum_{j=0}^{n-1}\partial_{\underline x}^{n-j-1}A_1\partial_{\underline x}^j\right),\quad n\ge2.
\end{equation}
We now proceed to examine the convergence of the series (\ref{series}) with the functions $A_n$ ($n\ge2$) given by (\ref{sol}). Let $\underline y$ be an arbitrary point in $\underline\Omega$. Then there exist a ball $B\big(\underline y,R(\underline y)\big)$ of radius $R(\underline y)$ centered at $\underline y$ and a positive constant $M(\underline y)$, such that
\[\left\vert\partial_{\underline x}^{n-j}A_s(\underline x)\partial_{\underline x}^{j}\right\vert\le M(\underline y)\frac{n!}{R^n(\underline y)},\;\;\underline x\in B\big(\underline y,R(\underline y)\big),\;\;j=0,\dots,n,\;\;s=0,1.\]
It follows that
\[\left\vert A_n(\underline x)\right\vert\le M(\underline y)\frac{n+R(\underline y)-1}{R^n(\underline y)},\quad\underline x\in B\big(\underline y,R(\underline y)\big),\]
and therefore the series (\ref{series}) converges normally in
\[\widetilde\Omega=\bigcup_{\underline y\in\underline\Omega}\left(-R(\underline y),R(\underline y)\right)\times B\big(\underline y,R(\underline y)\big).\]
Note that $\widetilde\Omega$ is a $x_0$-normal open neighbourhood of $\underline\Omega$ in $\mathbb R^{m+1}$, i.e. for each $x\in\widetilde\Omega$ the line segment $\{x+t:\,t\in\mathbb R\}\cap\widetilde\Omega$ is connected and contains one point in $\underline\Omega$.

We thus have proved the following.

\begin{theorem}
The function $\mathsf{CK}[A_0,A_1]$ given by
\begin{multline}\label{CK}
\mathsf{CK}[A_0,A_1](x)=A_0(\underline x)+x_0A_1(\underline x)\\
-\sum_{n=2}^\infty\frac{(-x_0)^n}{n!}\left(\sum_{j=0}^{n-2}\partial_{\underline x}^{n-j-1}A_0(\underline x)\partial_{\underline x}^{j+1}+\sum_{j=0}^{n-1}\partial_{\underline x}^{n-j-1}A_1(\underline x)\partial_{\underline x}^j\right)
\end{multline}
is inframonogenic in a $x_0$-normal open neighbourhood of $\,\underline\Omega$ in $\mathbb R^{m+1}$ and satisfies conditions (\ref{cond1})-(\ref{cond2}).
\end{theorem}

It is worth noting that if in particular $A_1(\underline x)=-\partial_{\underline x}A_0(\underline x)$, then
\[\mathsf{CK}[A_0,-\partial_{\underline x}A_0](x)=\sum_{n=0}^\infty\frac{(-x_0)^n}{n!}\,\partial_{\underline x}^nA_0(\underline x),\]
which is nothing else but the left monogenic extension (or CK-extension) of $A_0(\underline x)$. Similarly, it is easy to see that $\mathsf{CK}[A_0,-A_0\partial_{\underline x}](x)$ yields the right monogenic extension of $A_0(\underline x)$ (see \cite{BDS,DSS,D,S,SJ}).

Let $\mathsf{P}(k)$ ($k\in\mathbb N_0$ fixed) denote the set of all $\mathbb R_{0,m}$-valued homogeneous polynomials of degree $k$ in $\mathbb R^m$. Let us now take $A_0(\underline x)=P_k(\underline x)\in\mathsf{P}(k)$ and $A_1(\underline x)=P_{k-1}(\underline x)\in\mathsf{P}(k-1)$. Clearly,
\begin{multline*}
\mathsf{CK}[P_k,P_{k-1}](x)=P_k(\underline x)+x_0P_{k-1}(\underline x)\\
-\sum_{n=2}^k\frac{(-x_0)^n}{n!}\left(\sum_{j=0}^{n-2}\partial_{\underline x}^{n-j-1}P_k(\underline x)\partial_{\underline x}^{j+1}+\sum_{j=0}^{n-1}\partial_{\underline x}^{n-j-1}P_{k-1}(\underline x)\partial_{\underline x}^j\right),
\end{multline*}
since the other terms in the series (\ref{CK}) vanish. Moreover, we can also claim that $\mathsf{CK}[P_k,P_{k-1}](x)$ is a homogeneous inframonogenic polynomial of degree $k$ in $\mathbb R^{m+1}$.

Conversely, if $P_k(x)$ is a homogeneous inframonogenic polynomial of degree $k$ in $\mathbb R^{m+1}$, then $P_k(x)\vert_{x_0=0}\in\mathsf{P}(k)$, $\partial_{x_0}P_k(x)\vert_{x_0=0}\in\mathsf{P}(k-1)$ and obviously $\mathsf{CK}[P_k\vert_{x_0=0},\partial_{x_0}P_k\vert_{x_0=0}](x)=P_k(x)$.

Call $\mathsf{I}(k)$ the set of all homogeneous inframonogenic polynomials of degree $k$ in $\mathbb R^{m+1}$. Then $\mathsf{CK}[.,.]$ establishes a bijection between $\mathsf{P}(k)\times\mathsf{P}(k-1)$ and $\mathsf{I}(k)$.

It is easy to check that

\[P_k(\underline x)=P_k(\partial_{\underline u})\frac{\left\langle\underline x, \underline u\right\rangle^k}{k!},\quad P_k(\underline x)\in\mathsf{P}(k),\]
where $P_k(\partial_{\underline u})$ is the differential operator obtained by replacing in $P_k(\underline u)$ each variable $u_j$ by $\partial_{u_j}$. Therefore, in order to characterize $\mathsf{I}(k)$, it suffices to calculate $\mathsf{CK}\big[\left\langle\underline x, \underline u\right\rangle^ke_A,0\big]$ and $\mathsf{CK}\big[0,\left\langle\underline x, \underline u\right\rangle^{k-1}e_A\big]$ with $\underline u\in\mathbb R^m$.

A simple computation shows that

\begin{multline*}
\mathsf{CK}\big[\left\langle\underline x, \underline u\right\rangle^ke_A,0\big](x)=\left\langle\underline x, \underline u\right\rangle^ke_A\\
-\sum_{n=2}^k\binom{k}{n}(-x_0)^n\left\langle\underline x, \underline u\right\rangle^{k-n}\left(\sum_{j=0}^{n-2}\underline u^{n-j-1}e_A\underline u^{j+1}\right),
\end{multline*}

\begin{multline*}
\mathsf{CK}\big[0,\left\langle\underline x, \underline u\right\rangle^{k-1}e_A\big](x)=x_0\left\langle\underline x, \underline u\right\rangle^{k-1}e_A\\
-\frac{1}{k}\sum_{n=2}^k\binom{k}{n}(-x_0)^n\left\langle\underline x, \underline u\right\rangle^{k-n}\left(\sum_{j=0}^{n-1}\underline u^{n-j-1}e_A\underline u^j\right).
\end{multline*}

\subsection*{Acknowledgments}

The second author was supported by a Post-Doctoral Grant of \emph{Funda\c{c}\~ao para a Ci\^encia e a Tecnologia}, Portugal (grant number: SFRH/BPD/45260/2008).

\end{document}